\documentclass[10pt,conference]{IEEEtran}
\usepackage[left=0.62in,right=0.62in,top=0.75in]{geometry}
\usepackage{latexsym,todonotes}
\usepackage{graphicx}
\usepackage{float}
\usepackage{subfigure}
\usepackage{multirow}
\usepackage{amsmath, amssymb}
\usepackage{txfonts}
\usepackage{mathrsfs}
\usepackage{exscale}%
\usepackage{relsize}%
\usepackage{xcolor}
\usepackage{algorithm,algorithmicx}
\usepackage{algpseudocode} 
\usepackage{tikz,mathpazo,pgfplots}
\usepackage{cite}
\usepackage{flowchart}
\usetikzlibrary{arrows,calc,positioning}
\usepackage{tabularx}%
\usepackage{textcomp}
\usepackage{booktabs}%
\usepackage[switch]{lineno}%
\usepackage[bookmarks=false]{hyperref}
\usepackage{soul}
\pgfplotsset{compat=1.14}
\algdef{SE}[DOWHILE]{Do}{doWhile}{\algorithmicdo}[1]{\algorithmicwhile\ #1}%

\begin{document} 
\title{Optimal Trajectory Planning and Task Assignment for UAV-assisted Fog Computing\\
\thanks{This work is supported by  UIC Research Grant R201911 and Zhuhai Basic and Applied Basic Research Foundation Grant ZH22017003200018PWC.}
} 

\author{\IEEEauthorblockN{ Shuaijun Liu\IEEEauthorrefmark{1}, Jiaying Yin\IEEEauthorrefmark{2}, Zishu Zeng\IEEEauthorrefmark{3}, Jingjin Wu\IEEEauthorrefmark{1}\IEEEauthorrefmark{4}}
\IEEEauthorblockA{\IEEEauthorrefmark{1}Department of Statistics and Data Science, BNU-HKBU United International College, Zhuhai 519087, China \\
\IEEEauthorrefmark{2}Li Ka Shing Faculty of Medicine, The University of Hong Kong, Pokfulam, Hong Kong SAR, China \\
\IEEEauthorrefmark{3}School of Engineering, The Hong Kong University of Science and Technology, Kowloon, Hong Kong SAR, China \\
\IEEEauthorrefmark{4}Guangdong Provincial Key Laboratory of Interdisciplinary Research and Application for Data
Science \\
Email: p930005034@mail.uic.edu.cn; jyyin@connect.hku.hk; alan.zeng@connect.ust.hk; jj.wu@ieee.org}}

\maketitle

\begin{abstract}
Fog computing is an emerging distributed computing model for the Internet of Things (IoT). It extends computing and caching functions to the edge of wireless networks. Uncrewed Aerial Vehicles (UAVs) provide adequate support for fog computing. UAVs can not only act as a relay between mobile users and physically remote edge devices to avoid costly long-range wireless communications but also are equipped with computing facilities that can take over specific tasks. In this paper, we aim to optimize the energy efficiency of a fog computing system assisted by a single UAV by planning the trajectories of the UAV and assigning computing tasks to different devices, including the UAV itself. We propose two algorithms based on the classical Ant Colony and Particle Swarm Optimization techniques and solve the problem by continuous convex approximation. Unlike most existing studies where the trajectories are assumed to be straight lines, we account for the effect of obstacles, such as buildings, and deliberately avoid them during the trajectory planning phase. Through extensive simulation experiments, we demonstrate that our proposed approach can achieve significantly better energy efficiency than existing benchmark algorithms.

\begin{IEEEkeywords}
Fog Computing, Task Assignment, Unmanned Aerial Vehicles (UAV), Trajectory Planning, Optimization Algorithm
\end{IEEEkeywords}
\end{abstract}

\section{Introduction}


Driven by the concept of the Internet of Things (IoT)~\cite{Ashton1999}, a large quantity of data are produced and collected by terminal mobile devices such as smart phones and tablets. To improve the efficiency of executing relevant tasks that utilize these data, fog computing has been proposed to move the computational functions from data centers in a centralized location to the edge of the network~\cite{Shi2016}. Architectures with task computing capabilities, such as micro data centers~\cite{Greenberg2009} or base stations, have been deployed at the edge of the network that are geographically closer to mobile devices. By shortening the average distance of wireless transmissions, significant improvements in various aspects including energy efficiency, transmission reliability and latency control have been achieved since fog computing techniques are widely adopted. However, it is usually not commercially profitable to deploy edge devices in regions where the density of mobile devices (MDs) is relatively low. Such coverage holes may create difficulty in handling tasks generated by MDs in remote locations.

The deployment of unmanned aerial vehicles (UAVs) is considered as a promising addition for fog computing especially in the above-mentioned scenario~\cite{Hammouti2019}. The advantages of UAVs and fog computing are strong complements for each other. In a fog computing architecture, UAVs can be considered as special edge nodes that can act as both a wireless access point and a computing facility. Particularly, UAVs, with their mobility nature, offer flexible cloud-to-thing connectivity for MDs in different geographical locations~\cite{Mao2017}.

On the other hand, due to the hardware constraints, UAVs have limited battery capacities, and are generally less powerful in computing tasks compared to traditional edge devices. Therefore, it is crucial to plan different stages of the computing process, including assigning tasks to different devices, allocating power and channels for different transmissions, planning trajectories of UAVs, and determining the location for a UAV to transmit its targeting MD, in order to fully exploit the advantages of integrating UAVs in fog computing.

Optimization problems related to task assignment and trajectory planning in UAV-based mobile edge computing have received wide attention. 
For example, Zhao \textit{et al.} \cite{Zhao2022} proposed a  cooperative multi-agent deep reinforcement learning framework to jointly minimize the execution delays and energy consumption while considering trajectory design, computation task allocation, and limited communication resource.
The authors of \cite{Wei2019} solved a joint optimization problem incorporating both the energy consumption and task latency by three decision-making algorithms. Li \textit{et al.}~\cite{Li2020} planned UAV trajectories to minimize user transmission power given the service requirements by the Dinkelbach algorithm and successive convex approximation technique. 

For formulations that concurrently involve multiple dependent objectives such as latency control and energy saving, decomposing the joint optimization problem to a number of single-objective problems has been a particularly popular approach to reduce the computational complexity~\cite{Wei2022, Ning2021, Chen2021, Huang2022}. Another common assumption, that we will also follow in this paper, is that the UAV flies at the fixed altitude while taking care to avoid collision with obstacles such as buildings~\cite{Li2020, Wei2022, Ning2021, Wang2022}. 

This paper aims to propose a more practical optimization framework for the scenario of UAV-assisted computational offloading in fog computing, to fill in the above-mentioned gaps in existing research. The contributions of this paper are summarized as follows,
\begin{itemize}
    \item We consider a formulation that takes into account most practical issues in UAV-assisted fog computing, including 1) energy consumption for UAV movement, task computation and communication, 2) delay for task computation and communication, and 3) obstacle avoidance in the trajectory planning phase.
    \item We optimally allocate wireless channels for transmission between the UAV and MDs based on the number and sizes of tasks from each MD to improve the overall efficiency in transmission.
    \item We design an Ant Colony-based algorithm to plan the trajectory of the UAV. The algorithm will return an optimal trajectory that minimizes the energy consumption for the movement of the UAV while avoiding obstacles. 
    \item We propose a Particle Swarm Optimization based algorithm to solve the non-convex optimization problem that deals with task assignment, power allocation, and processing frequency distribution at a specific time instance. We then identify the optimal point on the trajectory for the UAV to stop moving and start transmission by a greedy approach that computes and compares relevant performance metrics if the UAV starts transmission at any point along the planned optimal trajectory. The outputs from the two algorithms will then be combined to obtain the final solution to the joint optimization problem.
    \item We demonstrate, by numerical simulation results, that our proposed method achieves significant improvement in energy efficiency compared to existing benchmark methods. Also, we show that our proposed method has a faster convergence speed than current state-of-art approaches.
    
\end{itemize}

The remainder of this paper is organized as follows. We will describe the model in Section \ref{sec: model} and our proposed algorithms in Section \ref{sec: algorithm}, respectively. The simulation design and numerical results are presented in Section \ref{sec: experiment}. Finally, Section \ref{sec: conclusion} concludes this paper.

\section{Network Model} \label{sec: model}
Let $\mathbb{R}_+$ and $\mathbb{N}_+$ represent the sets of positive reals and positive integers, respectively.

We consider a basic structure of a UAV-assisted fog computing network, where $K \in \mathbb{N}_+$ MDs (e.g., smartphones), a single UAV, and a remote data center are presented. The MDs are initiating computing tasks following a Poisson process. The tasks could be executed locally and offloaded to the UAV or the remote data center. The total number of available wireless channels for transmission among the UAV and MDs is $N_c\in\mathbb{N}_+$. 


\begin{figure}[!ht]
	\centering
	\includegraphics[width=0.43\textwidth]{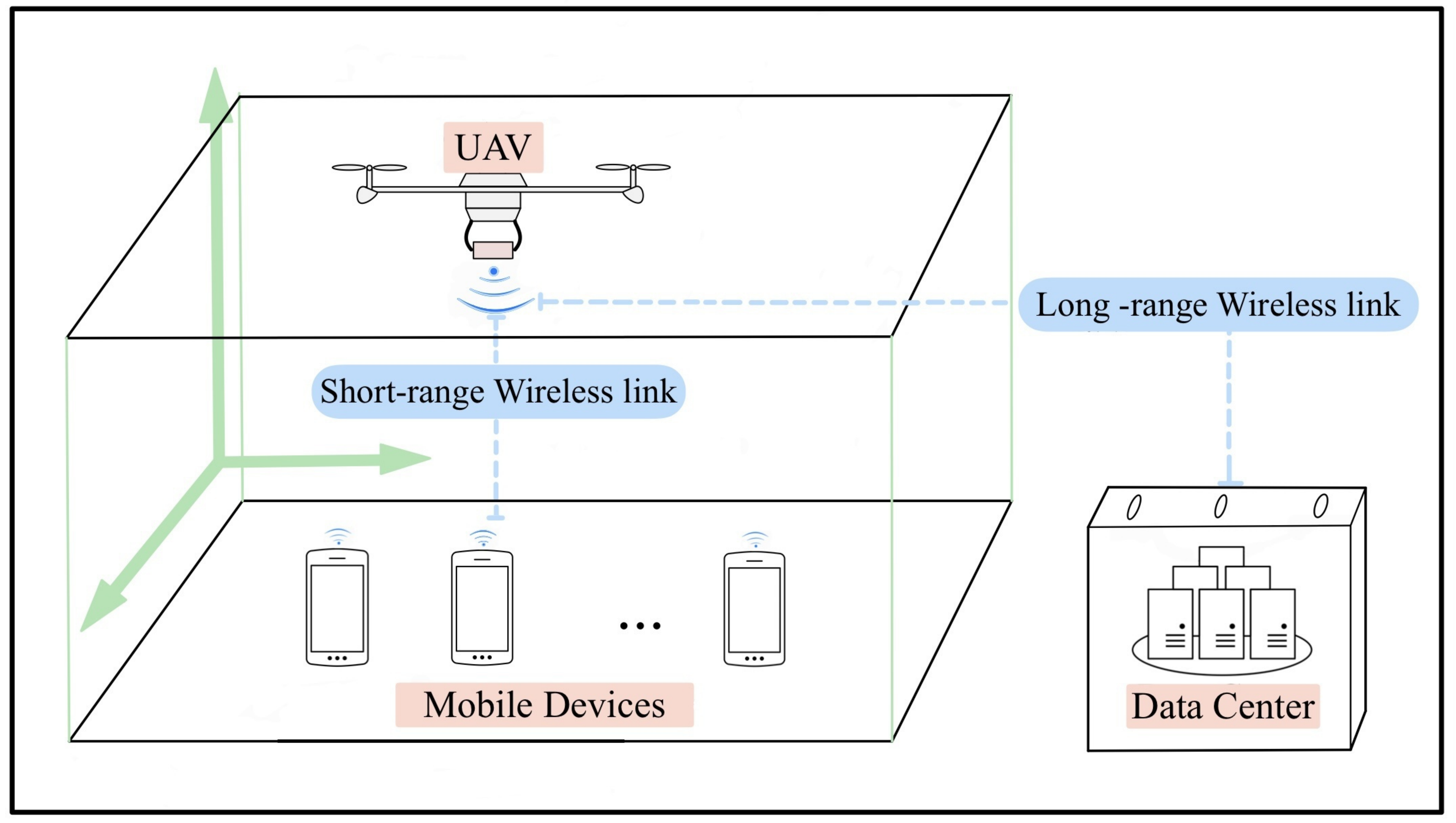}
	\caption{The structure of a UAV-assisted fog computing network.}
	\label{fig: UAV-fog structure}
\end{figure}

We illustrate the network structure in Fig. \ref{fig: UAV-fog structure}, where the UAV can establish short-range transmissions with MDs. Time is divided into slots, and the length of every timeslot is denoted as $L\in\mathbb{R}_+$. We assume that the UAV always flies at a constant altitude $h$, while all MDs are on the ground with an altitude of $0$. For ease of reference, we list the main notations used in this paper in Table \ref{tab: notation}.


\begin{table}[tbhp]
\vspace{-0.1cm}
\caption{Table of notations\label{tab: notation}}
\centering{}%
\renewcommand{\arraystretch}{1.5}
\begin{tabular}{l>{\raggedright}m{5.5cm}}
\toprule 
\textbf{Notation} & \textbf{Definition}\tabularnewline
\midrule
$K$ & The total number of MDs \tabularnewline
$N_c$ & The total number of wireless channels \tabularnewline
$N_j$ & The number of arrived tasks at the $j$th MD \tabularnewline
$N_o$ & The number of obstacles  \tabularnewline
$N_t$ & The number of turns taken by the UAV \tabularnewline
$L$ & The length of every timeslot \tabularnewline
$T$ & The total number of timeslots during the UAV to end all tasks and return\tabularnewline
$h$ & The flight altitude of the UAV \tabularnewline
$(x_j,y_j)$ & The location of the $j$th MD \tabularnewline
$(x(t),y(t))$ & The location of the UAV at the $t$th timeslot \tabularnewline
$C_j(t)$ & The number of channels assigned to the $j$th MD at the $t$th timeslot \tabularnewline
$v(t)$ & The velocity vector of the UAV at the $t$th timeslot \tabularnewline
$IP$ & The inertial factor of particles in the PSO \tabularnewline
$AP1,AP2$ & The acceleration factors of particles in the PSO \tabularnewline
$PH_i$ & The pheromones of ants in the ACO \tabularnewline
$HV_i$ & The heuristic values of ants in the ACO \tabularnewline
$\rho$ & The pheromone evaporation rate in the ACO \tabularnewline
$\mathcal{L}(t)$ & The weighted sum of network energy consumption and delay at the $t$th timeslot \tabularnewline
$p(t)$ & The transmission power of a Virtual Machine (VM) in the UAV at the $t$th timeslot \tabularnewline
$p^{max}(t)$,$p^{min}(t)$ 
 & Maximum and minimum transmission power of a VM in the UAV at the $t$th timeslot \tabularnewline
$p_j(t)$ & The $j$th MD's transmission power at the $t$th timeslot \tabularnewline
$p^{max}_{j}(t)$,$p^{max}_{j}(t)$ & Maximum and minimum transmission power for $j$th MD at the $t$th timeslot \tabularnewline
$f(t)$ & The processing frequency of a VM in the UAV at the $t$th timeslot \tabularnewline
$f^{max}(t)$,$f^{min}(t)$ & Maximum and processing frequency of a VM in the UAV at the $t$th timeslot \tabularnewline
$f_{j}(t)$ & The processing frequency of the $j$th MD at the $t$th timeslot \tabularnewline
$f^{max}_j(t)$,$f^{min}_j(t)$ & Maximum and minimum processing frequency of the $j$th MD at the $t$th timeslot \tabularnewline
\begin{tabular}[l]{@{}l@{}}$E^{\text{MD}}_{ij}(t)$, \\ $E^{\text{UAV}}_{ij}(t), E^{\text{DC}}_{ij}(t)$ \end{tabular} & The energy consumption of executing the $i$th task from the $j$th MD locally, at the UAV, or the DC at the $t$th timeslot \tabularnewline
\begin{tabular}[l]{@{}l@{}}$D^{\text{MD}}_{ij}(t)$, \\ $D^{\text{UAV}}_{ij}(t), D^{\text{DC}}_{ij}(t)$\end{tabular} & The delay of executing the $i$th task from the $j$th MD locally, at the UAV, or the DC at the $t$th timeslot \tabularnewline
\bottomrule

\end{tabular}
\end{table}

\subsection{UAV Movement Model} \label{subsec: UAV Movement model}
The communication between the UAV and an MD can be established only if they are sufficiently close to each other and there is a Line-of-Sight path between them. During the flight, the UAV should strictly comply with aviation control, and the flight altitude can not exceed the building. 
Based on the above requirements, we carry out optimal trajectory planning for the UAV currently completing the $i$th task. 
Therefore, the UAV must constantly move to attain different MDs' tasks at different times. The following equations can summarize the movement of the UAV.
Denote the location of the UAV at the $t$th timeslot as $(x(t), y(t))$ at the $t$th timeslot, and the flight altitude is a constant $h$\cite{Wang2022}. Given velocity vector $v(t)$, the travel distance of the UAV at the $t$th timeslot can be expressed as:

\begin{equation}
\vspace{-0.2cm}
\Delta(x(t), y(t))=\left(v_{x}(t), v_{y}(t)\right) \cdot L
\vspace{-0.05cm}
\end{equation}
\begin{equation}  
\vspace{-0.05cm}
d_{u}(t)= \sqrt{\left(x(t)-x(t-1)\right)^{2}+\left(y(t)-y(t-1)\right)^{2}}.
\end{equation}

It is desirable for the UAV to follow a shorter trajectory and make fewer turns due to power consumption and lifespan issues. In this regard, we define the trajectory planning value $\mathcal{R}$ as follows, 

\begin{equation}\label{eqn: path_planning}
\vspace{-0.05cm}
\mathcal{R}=\sum_{t=1}^{T} d_{u}(t) +\sum_{z=1}^{N_{t}}\left(\frac{\theta_{z}}{180}\right)^{\phi}, 
\end{equation}
where $T$ is the total number of timeslots during the UAV to end all tasks and return to charging, $N_t$ is the number of turns that the UAV makes, $\theta_z$ is the angle of $v_{u}$ before and after the $z$th turn, and $\phi$ is a coefficient determined by the acceleration at the time of turning. The value of $\mathcal{R}$ quantifies the requirements in trajectory planning. A smaller $\mathcal{R}$ will be achieved with moving distance, smaller turning angles, or fewer turns.

\subsection{Task Arrival and Channel Allocation Model} \label{subsec: Task Arrival and Channel Allocation Model}
Since the arrival pattern of tasks is not fully known, we build a task arrival model to realize the continuous auxiliary computing process that takes the dynamic arrival characteristics of the tasks into consideration. To ensure the randomness of tasks arrival, we assume that the task arrivals from $j$th MD conform to a Poisson process with an arrival rate of $\lambda_{j}$. For the $i$th task from the $j$th MD, we denote $c_{i j}$ as the number of CPU cycles required to process each input task, and $s_{i j}$, $o_{i j}$ as the input and output task data size, respectively. The values of $s_{i j}$ and $o_{i j}$ are random variables that follow the exponential distribution.

After determining the mission arrival mode, the key is to allocate the limited wireless channel. In order to allocate the wireless channels according to the task sizes efficiently, we assume that the proportion of channels allocated to the $j$th MD follows the Gamma probability distribution with shape parameter $\alpha$ and scale parameter $\beta$, denoted as,

\begin{equation}\label{eqn: channel}
\vspace{-0.05cm}
\omega_j(t)= \frac{\beta^{\alpha}}{\Gamma(\alpha)}s_{ij}^{\alpha-1}e^{\beta s_{ij}},
\end{equation}
where we set the parameters $\alpha=\beta=2$. The proportion of the probability density function value generated by taking the task sizes as the variable in the total quantity will be the proportion of the $\omega_{j}(t)$ in the $N_c$. Hence, the number of channels allocated to the $j$th MD at timeslot $t$ is $C_j(t) = N_c \cdot \omega_j(t)$.
The size of the task also called the storage capacity, will determine the number of channels required to transmit the task at the current moment. When the size of a certain arrival task $s_{i j}$ is too large, it will occupy too many computing resources and reduce the execution efficiency of UAV. Such a setting provides more opportunities for the task with great demand to transmit channels. At the same time, it is punitive to a single task with a large storage capacity to prevent the channel from being monopolized and guarantee a certain level of transmission efficiency.



\subsection{Task Assignment Optimization Model} \label{subsec: Task Assignment Optimization Model}
At the $t$th timeslot, for tasks from the $j$th MD, we define the utility function considering delay and energy consumption for transmission and computation as,

\begin{equation*}
\centering
\begin{split}
    \mathcal{D}_{j}(t) = \sum_{i=1}^{N_j}\left((D_{ij}^{\text{MD}}(t) \cdot x_{ij}^{\text{MD}}(t) + D_{ij}^{\text{UAV}}(t) \cdot x_{ij}^{\text{UAV}}(t) + \right. \\
    D_{ij}^{\text{DC}}(t) \cdot x_{ij}^{\text{DC}}(t)) + \epsilon (E_{ij}^{\text{MD}}(t) \cdot x_{ij}^{\text{MD}}(t)  \\
\left. + E_{ij}^{\text{UAV}}(t) \cdot x_{ij}^{\text{UAV}}(t) + E_{ij}^{\text{DC}}(t) \cdot x_{ij}^{\text{DC}}(t))\right),
\end{split}
\label{eqn:obj_function}
\end{equation*}
where we define $x_{ij}(t)=(x_{ij}^{\text{MD}}(t), x_{ij}^{\text{UAV}}(t), x_{ij}^{\text{DC}}(t))$ as an array of binary variables indicating the executing location of the $i$th task from the $j$th MD at the $t$th timeslot. $\epsilon$ is a weighting factor that accounts for the relative importance between energy consumption and delay in the utility function. 
When the tasks are more urgent in nature and thus require a shorter delay, $\epsilon$ may take a relatively small value to emphasize the contribution of delay to the utility function. Conversely, when the energy is short in supply, and the main objective is to reduce energy consumption, we can set a higher $\epsilon$. 
Then, we define $\mathcal{S}(t) = \sum_{j=1}^{K} \left( \mathcal{D}_{j}(t)\right)$  as the sum of the utility function for all MDs.

We aim to minimize the total consumption of the UAV-assisted fog computing network, including the consumption for completing the tasks and the consumption for the UAV movement. In other words, we also need to account for the energy consumed by the UAV for its movement. The joint optimization problem is formulated as,  
\begin{equation}
\vspace{-0.1cm}
    \begin{array}{ll}
\underset{ x_{ij}}{\operatorname{Minimize}} & \mathcal{L}(t) = \mathcal{S}(t) + \epsilon E_{u}(t) \\
\text { subject to } 
& v(t) \leq v^{\max} \\
& \sum_{j=1}^{K} C_j(t) \leq N_{c}\\
& f^{\min } \leq f(t) \leq f^{\max } \\
& f_{j}^{\min } \leq f_{j}(t) \leq f_{j}^{\max } \\
& p^{\min } \leq p(t) \leq p^{\max } \\
& p_{j}^{\min } \leq p_{j}(t) \leq p_{j}^{\max } \\
& x_{i j}^{\text{MD}}(t)+x_{ij}^{\text{UAV}}(t)+x_{ij}^{\text{DC}}(t)=1 \\
& x_{i j}^{\text{MD}}(t), x_{ij}^{\text{UAV}}(t), x_{ij}^{\text{DC}}(t) \in\{0,1\} \\
& j \in\{1, \ldots, K\},
\end{array}\label{eq:problem}
\end{equation}
where $E_{u}(t)$ = 0.5$M \cdot L \cdot \Vert v(t)\Vert ^{2} $ is the energy consumption for the movement of the UAV, and $M$ is the mass of the UAV. Note that processing frequencies $f(t)$, $f_{j}(t)$, and transmission powers $p(t)$, $p_j(t)$ are all functions of the decision variables $x_{ij}(t)$. The frequencies and powers in turns would determine the total consumption. Detailed relationships can be found in a number of existing literature, e.g.,~\cite{Wei2019}. 

Note that Problem~\eqref{eq:problem} is by nature non-convex. Therefore, the traditional approaches to solve it are computationally prohibitive.

\section{Algorithms}
\label{sec: algorithm}
\subsection{Overview}
    Our proposed solution to Problem~\eqref{eq:problem} mainly consists of two global optimization algorithms. The first algorithm focuses on assigning tasks and allocating resources, including transmission power and frequency, by the Particle Swarm Optimization technique. In contrast, the second determines the UAV's optimal trajectory based on Ant Colony Optimization. Finally, the optimal trajectory and resource allocations at every point on the trajectory, namely the output of the first two algorithms, are integrated to determine the optimal position for the UAV to start transmission with the MD.


\subsection{Task Assignment, Power and Frequency Allocations}

The $i$th task from the $j$th MD at the $t$th timeslot can be executed in the $j$th MD itself locally, the UAV, or the remote data center. To ensure that all the computational resources are efficiently utilized during the transmission and computation of user tasks, it is essential to deal with the task scheduling policy and the corresponding transmission parameter settings. 

We consider solving this problem by Particle Swarm Optimization\cite{Zhang2015}. As mentioned before, assignment decisions are represented by $x_{ij}=(x_{ij}^{\text{MD}}, x_{ij}^{\text{UAV}}, x_{ij}^{\text{DC}})$. For the convenience of presentation, we group the task assignment decision, transmission power and processing frequency together as a particle group $s_m = (x_{ij}, p_j(t), f_j(t))$ in the particle swarm, where $m \in \{1, ..., M\}$. Then, we decide the location of the particle groups by uniformly sampling $M\in\mathbb{N_+}$ particle groups. To find the optimal solution iteratively, we initialize the velocity of the particle groups and update according to the method in \cite{Zhang2015}. We terminate the algorithm when the difference between results in two consecutive iterations is smaller than a threshold $\xi$. The algorithm outputs the minimum value $\mathcal{S}^{*}$, the optimal task assignment decision $x_{ij}^{*}$, transmission power $p_j^*(t)$, and processing frequency $f_j^*(t)$.



\renewcommand{\algorithmicrequire}{\textbf{Input:}}  
\renewcommand{\algorithmicensure}{\textbf{Output:}}
\begin{algorithm}[ht]
	\caption{(PSO) Task Assignment, Power and Frequency Allocations based on Particle Swarm Optimization} 
	\begin{algorithmic}[1]
		\Require The intervals of $x_{ij}$, $p_j(t)$, $f_j(t)$, $AP1$, $AP2$, $IP$
		\Ensure $\mathcal{S}^{*}(t)$; $x_{ij}^{*}$; $p_j^{*}(t)$; $f_j^{*}(t)$
		\vspace{1 ex}
		\For {$m = 1$ to $M$}
			\State Initialize partilces' location  $X_m$
			\State Initialize velocity $V_m$ 
			\State Initialize $pBest_m$ to its location $pBest_m \longleftarrow X_m$
		\EndFor
	    \State Initialize $gBest(0)$ $\longleftarrow$ argmin$fit(pBest_m)$, where $fit$ represent the equation  to compute $\mathcal{S}^{*}(t)$; $n = 1$
	   \Do
	    \For{$m=1$ to $M$}
			\State Update $V_m$ and $X_m$ by acceleration factors 
			\Statex \quad \quad \quad $AP1$, $AP2$ and $IP$
			\If{$fit(X_m) < fit(pBest_m)$}
 			   \State $pBest_m \longleftarrow X_m$
	            \If{$fit(pBest_m)<fit(gBest)$}
 		        \State $gBest(n) \longleftarrow pBest_m$
 		   	    \EndIf
	        \EndIf
	        \State $n \longleftarrow (n+1)$
		\EndFor
	   \doWhile{$\lvert gBest(n+1)-gBest(n) \rvert < \xi$}
	    \State $\mathcal{S}^{*}(t) \longleftarrow fit(gBest)$
	    \State $(x_{ij}^{*}$; $p_j^{*}(t)$; $f_j^{*}(t))$ $\longleftarrow gBest$
	    \State Output $\mathcal{S}^{*}(t)$; $(x_{ij}^{*}$; $p_j^{*}(t)$; $f_j^{*}(t))$
	\end{algorithmic} \label{Algorithm1}
	\vspace{-0.1cm}
\end{algorithm}

\subsection{UAV Optimal Trajectory Planning}
Assume that there are $N_o \in [N_o^{\text{min}}, N_o^{\text{max}}]$ obstacles randomly deployed within the region. We use $\boldsymbol O = \{(x^o_i, y^o_i)\}$ to denote the set of locations of all obstacles, where $i \in \{1, 2, ..., N_o\}$. The heights of all the obstacles are larger than or equal to $h$, so the UAV needs to avoid all of them. Avoiding obstacles can be considered a constraint for the UAV in the optimization problem, and changing the direction at the $t$th timeslot could ensure that the coordinates of the UAV in the next timeslot will not fall within the coordinates of the obstacles. 

We compute the movement trajectory by Algorithm~\ref{Algorithm2} based on the Ant Colony Algorithm\cite{Dorigo2006} to minimize $\mathcal{R}$, so as to get the optimal movement path. We first initialize the pheromones $PH_0$ and heuristic values $HV$ of the ants, and the evaporation rate $\rho$. While implementing the algorithm, the UAV could be considered as an ant, and we record the motion state and the location of the ant at each timeslot, the number of UAV turnings $\phi$, as well as the turning angle $\theta_z$, velocity $v_u(t)$, and acceleration at each turn. Finally, the minimum trajectory planning value $\mathcal{R}$ is obtained through several iterations, and the information of each coordinate through which the ant moves constitutes the final output, namely the optimal trajectory of the UAV.

\renewcommand{\algorithmicrequire}{\textbf{Input:}}  
\renewcommand{\algorithmicensure}{\textbf{Output:}} 
\begin{algorithm}[ht]
\vspace{-0.05cm}
	\caption{(ACO) Trajectory Planning based on Ant Colony Optimization} 
	\begin{algorithmic}[1]
		\Require
		$(x_j,y_j)$; $[N_o^{\min}$,$N_o^{\max}]$; $\rho$; $HV$
		\Ensure
		$\mathcal{R}_{min}$; $(x_k(t),y_k(t))$ 
		\vspace{1 ex}
		\State Randomly generate $N_o \in [N_o^{\min}$,$N_o^{\max}]$ obstacles
		\Statex record the location as coordinates $(x^o_i, y^o_i) \in \boldsymbol O$
        \For {$m = 1$ to $M$}
        \Do
            \State Randomly set initial coordinate $(x(0),y(0))$
        \doWhile{$(x(0),y(0)) \in \boldsymbol O$}
		\For{each edge}
			\State Set initial pheromone $PH$
		\EndFor
		\For{each ant $k$}
			\State Initial coordinates $(x_k(t),y_k(t))$ = $(x(0),y(0))$
	    \For{each edge}
	        \Do
	            \State Choose the next coordinate with the 
	            \Statex \quad \quad \quad \quad \quad \quad probability by $PH$ and $HV$
	        \doWhile{$(x_k(t+1),y_k(t+1)) \in \boldsymbol O$}
	    	\State output $(x_k(t+1),y_k(t+1))$
		\EndFor
		\State Compute and output the length $\sum_{t=1}^{T} d_{u}(t)$ of    
		\Statex \quad \quad \quad the path by the $k$th ant and the $\mathcal{R}$ value
		\For {each edge}
	       \State Update the $i$th pheromone value $PH$ by $\rho$
        \EndFor
        \EndFor
        \EndFor
		\State Compute and output the $\mathcal{R}_{min}$ by equation \eqref{eqn: path_planning}.
	\end{algorithmic} \label{Algorithm2}
\end{algorithm}
\subsection{Optimal Transmission Position of the UAV}
With the optimal trajectory determined previously, we can confirm the coordinates of the UAV at each timeslot, assuming that it sticks to the optimal trajectory. A reasonable assumption is that the UAV will remain at a fixed location once it starts transmission with the targeted MD, in order to guarantee the reliability of the transmission. 

We consider a straightforward approach to determine the optimal position on the optimal trajectory for the UAV to stop moving and start transmission. We iteratively check the coordinate ($x^*(t)$, $y^*(t)$) for every timeslot $t = 1,2,\cdots$ on the optimal trajectory obtained by Algorithm \ref{Algorithm2}, and invoke Algorithm~\ref{Algorithm1} to compute the consumption $\mathcal{S}^{*}(t)$ with the optimal task assignment and resource allocation if the UAV is transmitting at ($x^*(t)$, $y^*(t)$). Finally, we compare the values of $\mathcal{L}^{*}(t)$ for all $t = 1,2,\cdots$, and identify the optimal timeslot $t^*$ for the UAV to start transmission as $t^* = \arg\min_t \mathcal{L}^{*}(t)$.

The procedures to obtain the optimal timeslot and corresponding coordinates for the UAV to start transmission are summarized in Algorithm~\ref{Algorithm3}.

\renewcommand{\algorithmicrequire}{\textbf{Input:}}  
\renewcommand{\algorithmicensure}{\textbf{Output:}}
\begin{algorithm}[!ht]
	\caption{(TDO) Travel Distance Optimization } 
	\begin{algorithmic}[1]
	\Require
		$\{(x^*(t),y^*(t))|t = 1,2,\cdots, T\}$
		\Ensure
		$\mathcal{L^*}(t)$; $t^*$; $(x^*(t^*),y^*(t^*))$
		\vspace{1 ex}
	    \For{each $t$}
	    \State Retrieve $(x(t),y(t))$ from  Algorithm~\ref{Algorithm2}
	    \State Invoke Algorithm~\ref{Algorithm1} to compute $\mathcal{S}^{*}(t)$ based on $(x^*(t),y^*(t))$
	    \State Compute $\mathcal{L}^{*}(t)$ based on $(x(t),y(t))$
	    \EndFor
	    \State Obtain $t^* = \arg\min_t \mathcal{L}^{*}(t)$
	    \State \Return $\mathcal{L}^{*}(t)$; $t^*$; $x^*(t^*),y*(t^*)$
	\end{algorithmic} \label{Algorithm3}
\end{algorithm}


%

\section{Performance Evaluation} \label{sec: experiment}
\subsection{Experiment setup}
In this section, we perform numerical simulations on systems with a range of parameter values to evaluate the effectiveness and adaptability of our proposed solutions. We consider that $K$ MDs are deployed in an area of $S \times S$. For simplicity without loss of generality, we discretize the area into $1 \times 1$ grids. We divide the time into multiple timeslots of equal length. The starting horizontal and vertical coordinates of the UAV are independently and randomly generated in [0, $S/2$]. Values of system parameters in the experiment are listed in Table \ref{tab:MODEL PARAMETER SETTINGS}. 

\begin{table}[htbp]
\caption{SYSTEM PARAMETER SETTINGS\label{tab:MODEL PARAMETER SETTINGS}}
\begin{center}
\begin{tabular}{c c c c}
\hline Parameter & Value & Parameter & Value\\
\hline $S$ & 10000$\mathrm{~m}$ & $\epsilon$ & [0.05, 1.00] \\
\hline
$h$ & 50$\mathrm{~m}$ & $C_j(t)$ & [0.05, 2.50] \\
\hline
$K$ & 50 & [$p^{\min}$,$p^{\max}$] & [40, 80]$\mathrm{~mW}$ \\
\hline
$N_{c}$ & 40 & [$p_{j}^{\min}$, $p_{j}^{\max }$] & [30, 70]$\mathrm{~mW}$ \\
\hline
$L$  & 0.1$\mathrm{~s}$ 
& [$f^{\min}$, $f^{\max}$] & [1.0, 2.0]$\mathrm{~GHz}$ \\
\hline
$v^{\max}$ & 10$\mathrm{~m/s}$ 
& [$f_{j}^{\min}$, $f_{j}^{\max }$] & [0.5, 2.0]$\mathrm{~GHz}$ \\
\hline
$\xi$ & 0.01 
& [$N_o^{\min}$, $N_o^{\max}$] & [2000, 3000] \\
\hline
\end{tabular}
\end{center}
\vspace{-0.5cm}
\end{table}

Based on the scenarios above and model parameter settings in Table~\ref{tab:MODEL PARAMETER SETTINGS}, we perform each algorithm to solve the joint optimization problem 50 times. The following results presented in this section are based on the average of the 50 runs for each corresponding method. In each run, the weighting factor $\epsilon$ is generated randomly within its domain.
The values of relevant parameters involved in the algorithms are listed in Table \ref{tab:ALGORITHM PARAMETER SETTINGS}. 

\begin{table}[htbp]
\caption{ALGORITHM PARAMETER SETTINGS \label{tab:ALGORITHM PARAMETER SETTINGS}}
\begin{center}
\begin{tabular}{cc|cc}
\hline Parameter & Value & Parameter & Value\\
\hline
$\rho$ & 0.25
& $PH_0$ & 3.8 \\
$AP1$ & 2.0
& $IP$ & 0.65 \\
$AP2$ & 2.0
& $HV$ & 2.5 \\
\hline
\end{tabular}
\end{center}
\vspace{-0.5cm}
\end{table}
We demonstrate and compare the results in five different scenarios with the following approaches. Firstly, we consider an obstacle-free area where the UAV always flies in the direction of the line connecting the origin and the destination. Next, we focus on the performance improvement achieved by optimally allocating the task computing position, allocating wireless channels, and determining the processing frequency and transmission power of MDs and UAV. We will present the results from the following scenarios with corresponding optimization approaches:

\begin{itemize}
\item RAN: The channels are all randomly allocated to $K$ MDs. All tasks are randomly assigned to the local MD, the UAV, or the data center. Processing frequencies and power allocations are generated randomly in the domains. 

\item GA: The Genetic Algorithm used in~\cite{Wei2019} is deployed to obtain task assignments, as well as the processing frequencies and transmission powers of MDs and the UAV. The other parameters are still randomly generated in the domains.

\item PSO: We apply Algorithm \ref{Algorithm1} (Particle Swarm Optimization) to determine task assignments, the processing frequencies ,and transmission powers of MDs and the UAV. The other parameters are still randomly generated in the respective ranges.

\item CA: On top of PSO, The channels are allocated according to the arrival rates a d sizes of tasks generated from different MDs, according to the method described in Section~\ref{subsec: Task Arrival and Channel Allocation Model}.

\item TDO: On top of CA, Algorithm \ref{Algorithm3} (Travel Distance Optimization) is deployed to determine the location for the UAV to transmit with the targeting MD. 

\end{itemize}

\subsection{Numerical results}

\subsubsection{Convergence of Algorithms}

\begin{figure}[ht]
    \vspace{-0.2cm}
    \centering
    \includegraphics[width=8.8cm]{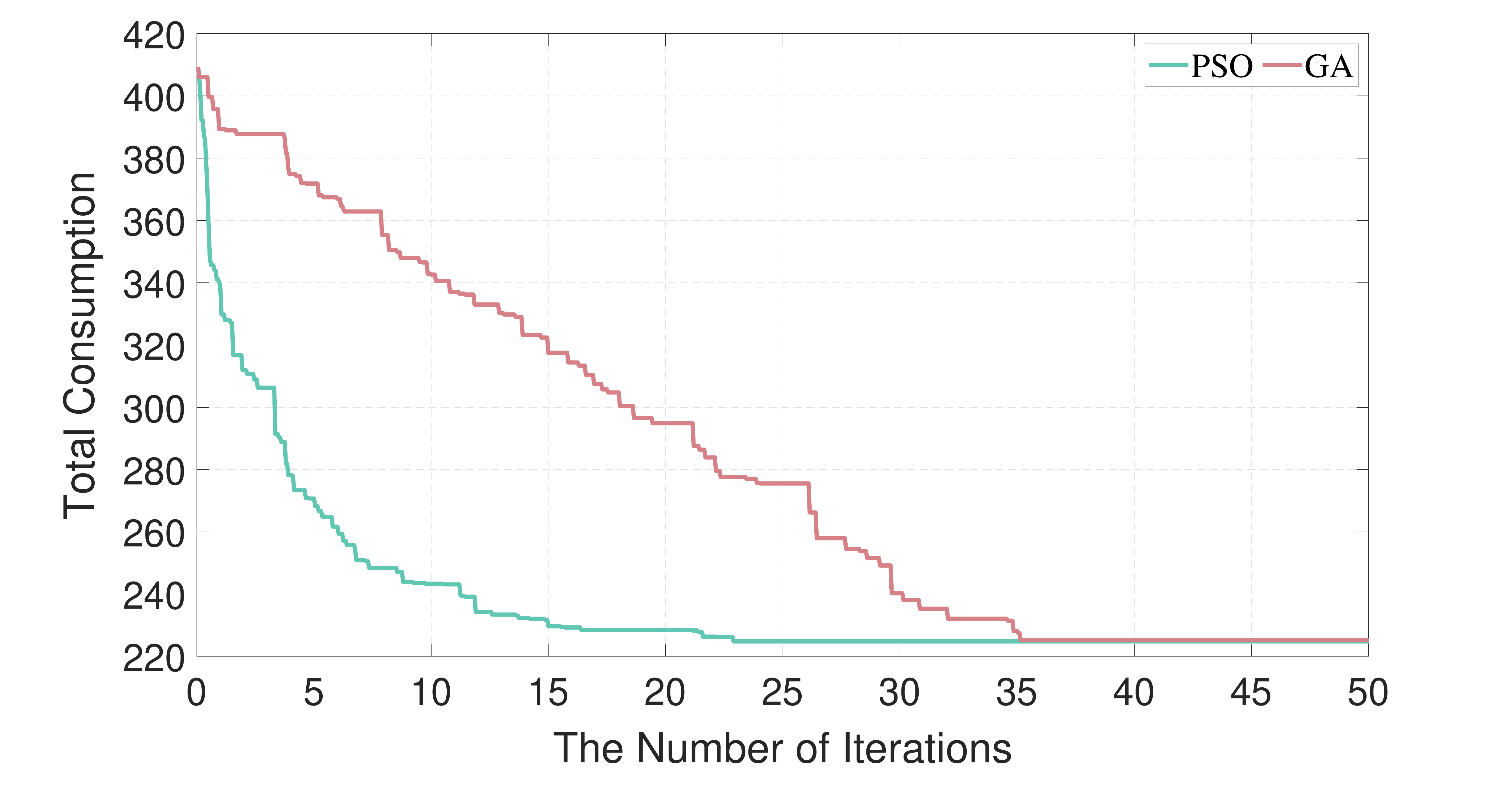}
    \caption{Convergence curves of PSO and GA.} 
    \label{fig: Iterations}
    \vspace{-0.2cm}
\end{figure}

We first verify our proposed approach's convergence and compare its convergence rate to the benchmark GA. The convergence curves of PSO and GA are shown in Fig.~\ref{fig: Iterations}, where the horizontal axis represents the number of iterations,
and the vertical axis denotes the total consumption of respective algorithms at a specific iteration. As shown by the curves, while the convergent values of the two algorithms are extremely close, PSO can achieve convergence in about $25$ iterations compared to $35$ iterations for GA. This represents an improvement of $28.6\%$ in convergence speed. 

\subsubsection{Task Assignments}

\begin{figure}[ht]
    \vspace{-0.2cm}
    \centering
    \includegraphics[width=8.8cm]{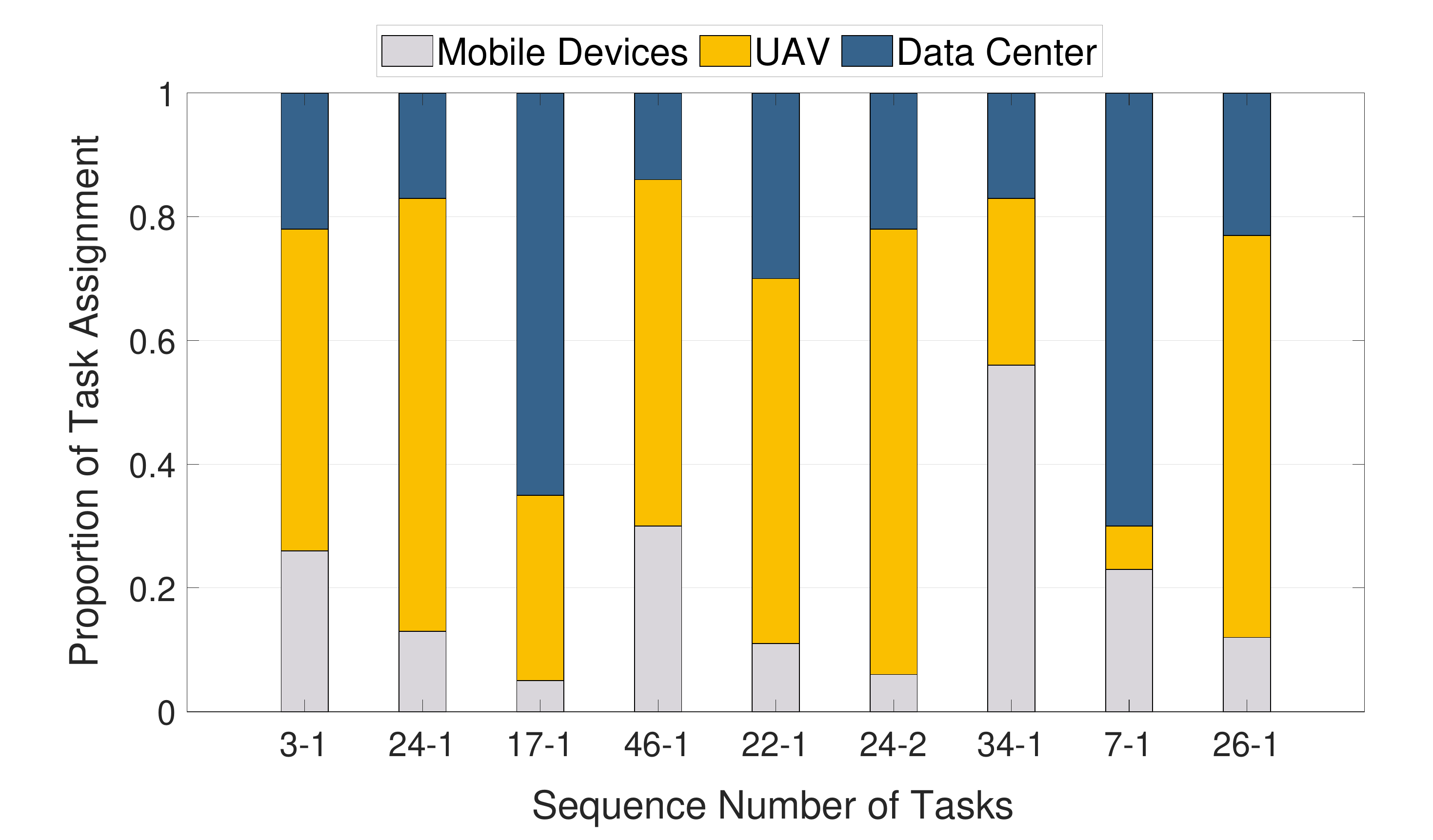}
    \caption{Assignments for $9$ tasks at one particular timeslot.} 
    \label{fig: TaskAssign}
    \vspace{-0.2cm}
\end{figure}

We consider a sample of $9$ tasks that are initiated at the same timeslot and illustrate their respective assignments by Algorithm~\ref{Algorithm1} in  Fig.~\ref{fig: TaskAssign}. Here, we denote the $i$th task from the $j$ MD as the $j$-$i$ task. For example, ``$3$-$1$" denotes the $1$st task from the $3$rd MD. Different colors in each column represent the proportion of a certain task that is executed in the MD, the UAV, and the cloud data center, respectively. 



\subsubsection{Channel Assignment}

\begin{figure}[ht]
    \vspace{-0.2cm}
    \centering
    \includegraphics[width=8.8cm]{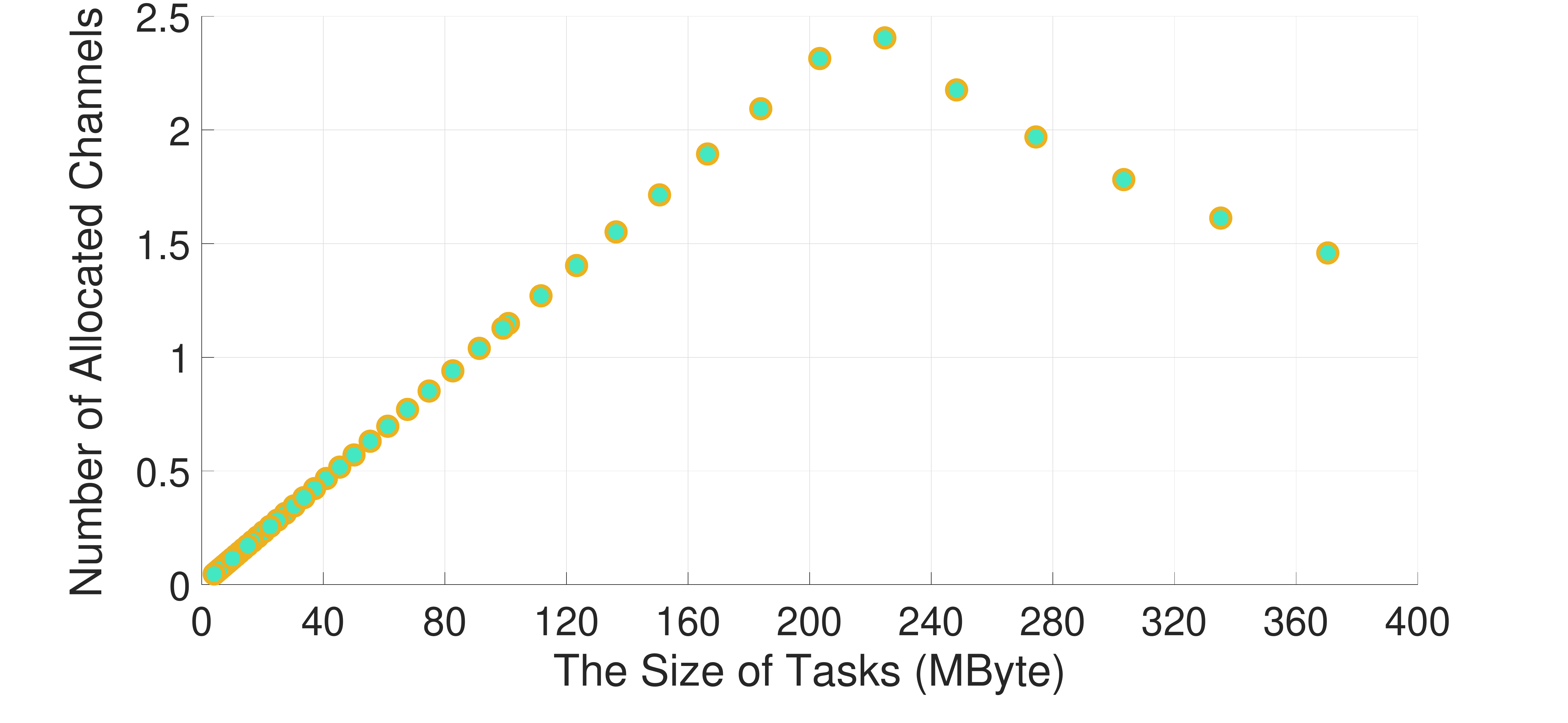}
    \caption{Channel allocation strategy for $36$ tasks of different sizes.} 
    \label{fig: Channel}
    \vspace{-0.2cm}
\end{figure}

Fig. \ref{fig: Channel} demonstrates our proposed channel assignment results with 36 tasks of different sizes at a time slot. Each point in the Fig. ~\ref{fig: Channel} indicates the number of channels (vertical axis) allocated to a task of a certain size (horizontal axis). The savings achieved by our proposed channel allocation strategy in total consumption is also shown by the curve CA compared to PSO in Fig.~\ref{fig: result1}.  

\subsubsection{Total Consumption}
\begin{figure}[ht]
    \vspace{-0.2cm}
    \centering
    \includegraphics[width=8.8cm]{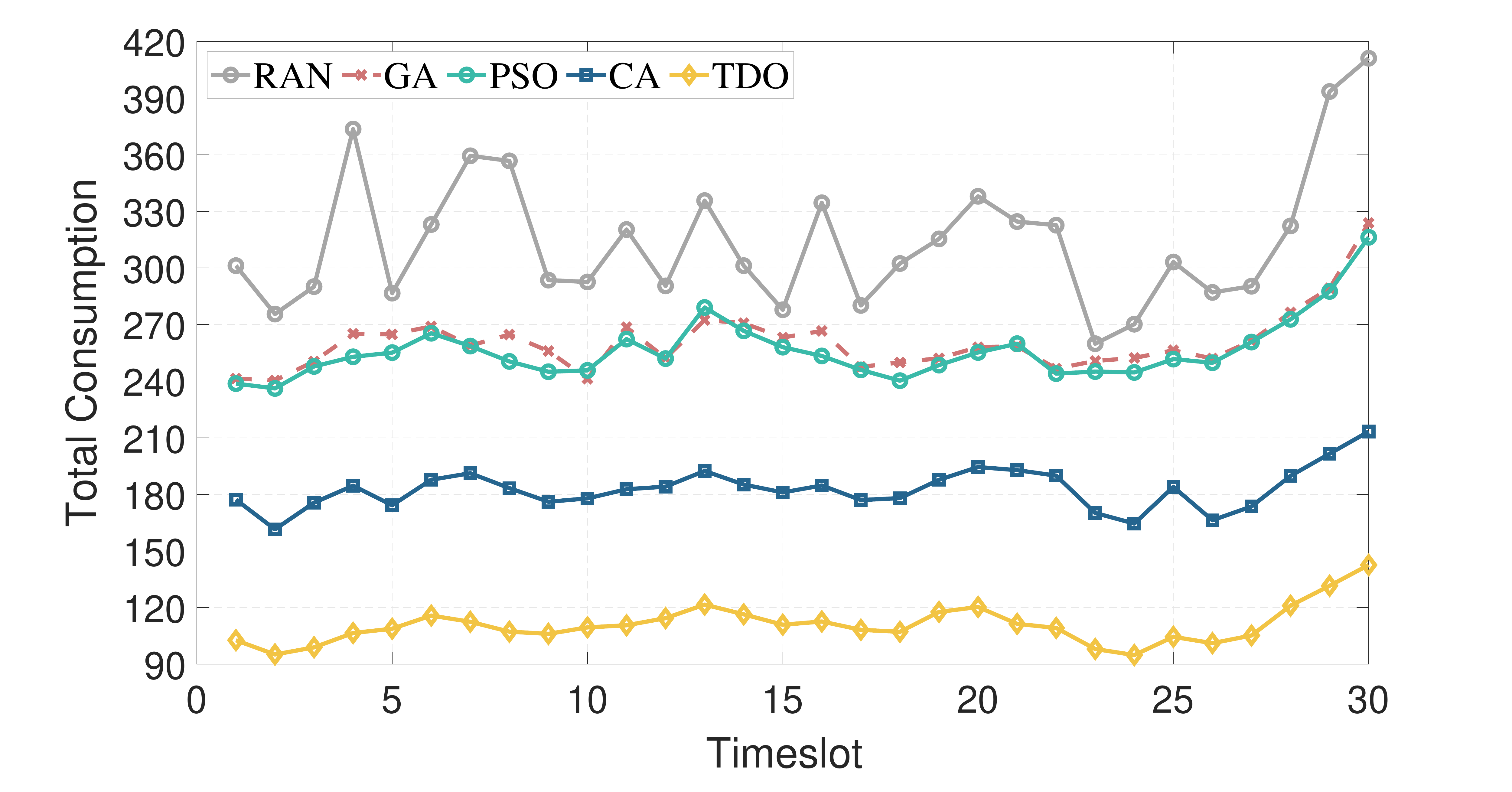}
    \caption{Comparison of total consumption by different methods} 
    \label{fig: result1}
    \vspace{-0.2cm}
\end{figure}

The results of total consumption achieved by the five approaches introduced in the previous subsection at different timeslots are presented in Fig. \ref{fig: result1}. We can see that GA and PSO's performances are very close to each other and much better than the baseline case (RAN). An optimized channel allocation by CA can further reduce the total consumption by up to $30\%$ compared to GA and PSO. The best of the five is the case where TDO optimizes the UAV transmission location on top of CA, reducing the total consumption by more than $45\%$.

\subsubsection{Trajectory planning}

\begin{figure}[htp!]
    \vspace{-0.2cm}
    \centering
    \includegraphics[width=7.6cm]{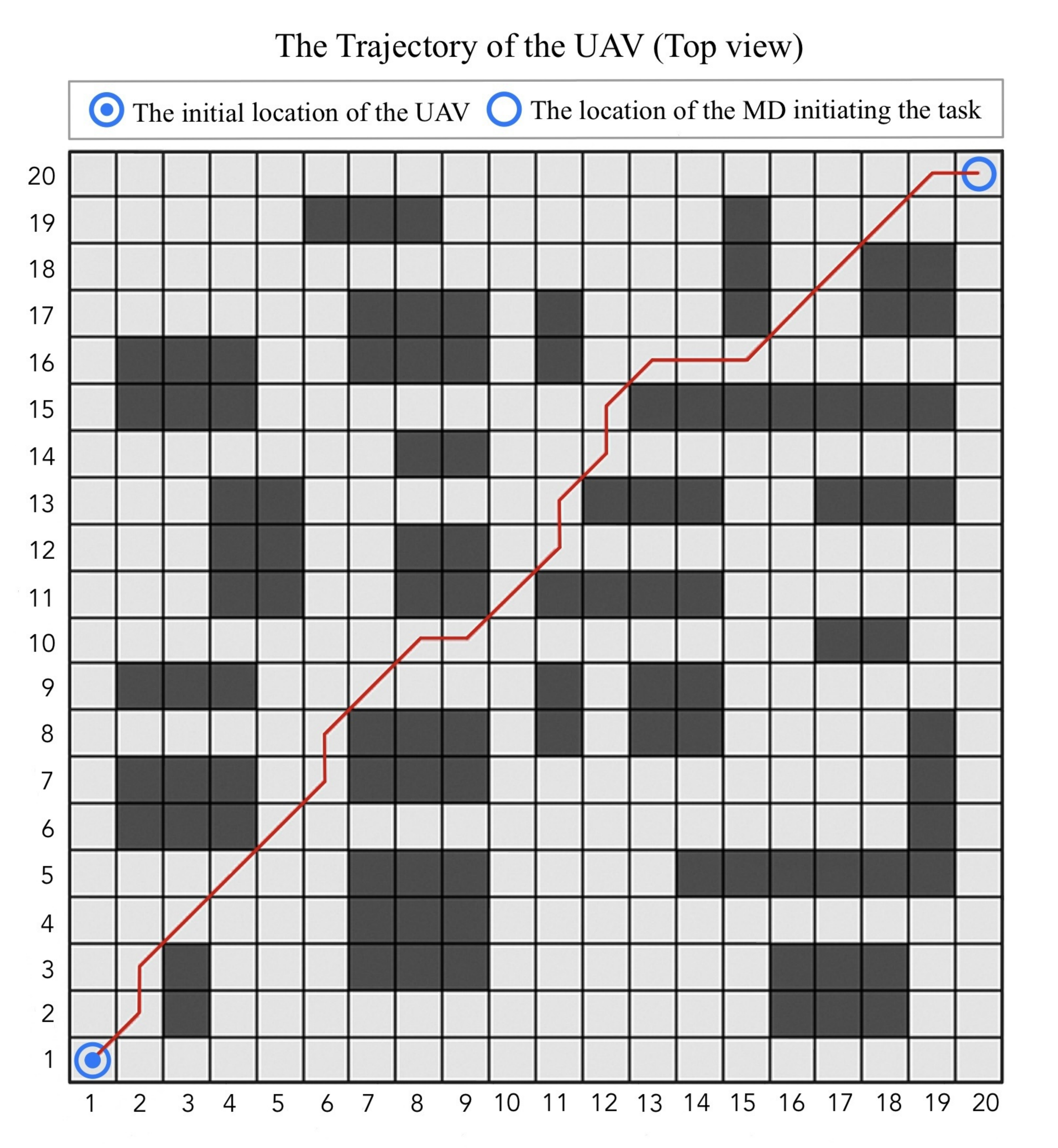}
    \caption{A demonstration of trajectory planning in an area with obstacles} 
    \label{fig: UAVpath}
\end{figure}

We then consider the effect of obstacles (such as office buildings) and test the performance of our trajectory planning algorithm.
In the given area, $N_o \in [N_o^{min}, N_o^{max}]$ obstacles are now randomly distributed, and every obstacle occupies one entire grid. After $200$ iterations, our algorithm obtains the optimal trajectory of the UAV movement. As shown in Fig. \ref{fig: UAVpath}, the dark grids represent obstacles that the UAV cannot fly over, and the light grids represent free space.

As demonstrated earlier, the trajectory planning algorithm will identify the optimal path between the starting location of the UAV and the destination where the MD initiating the task is located while avoiding the obstacles. Fig. \ref{fig: UAVpath} also shows an example of the planned trajectory. The side length of each square represents one unit of distance, the black square is the obstacle, and the red line is the best trajectory of the UAV. Note that according to our settings, the UAV does not necessarily need to travel the entire trajectory and reach the exact location of the MD, 
it may stop moving halfway and start transmission immediately.

\begin{figure}[htbp]
    \vspace{-0.2cm}
    \centering
    \includegraphics[width=8.8cm]{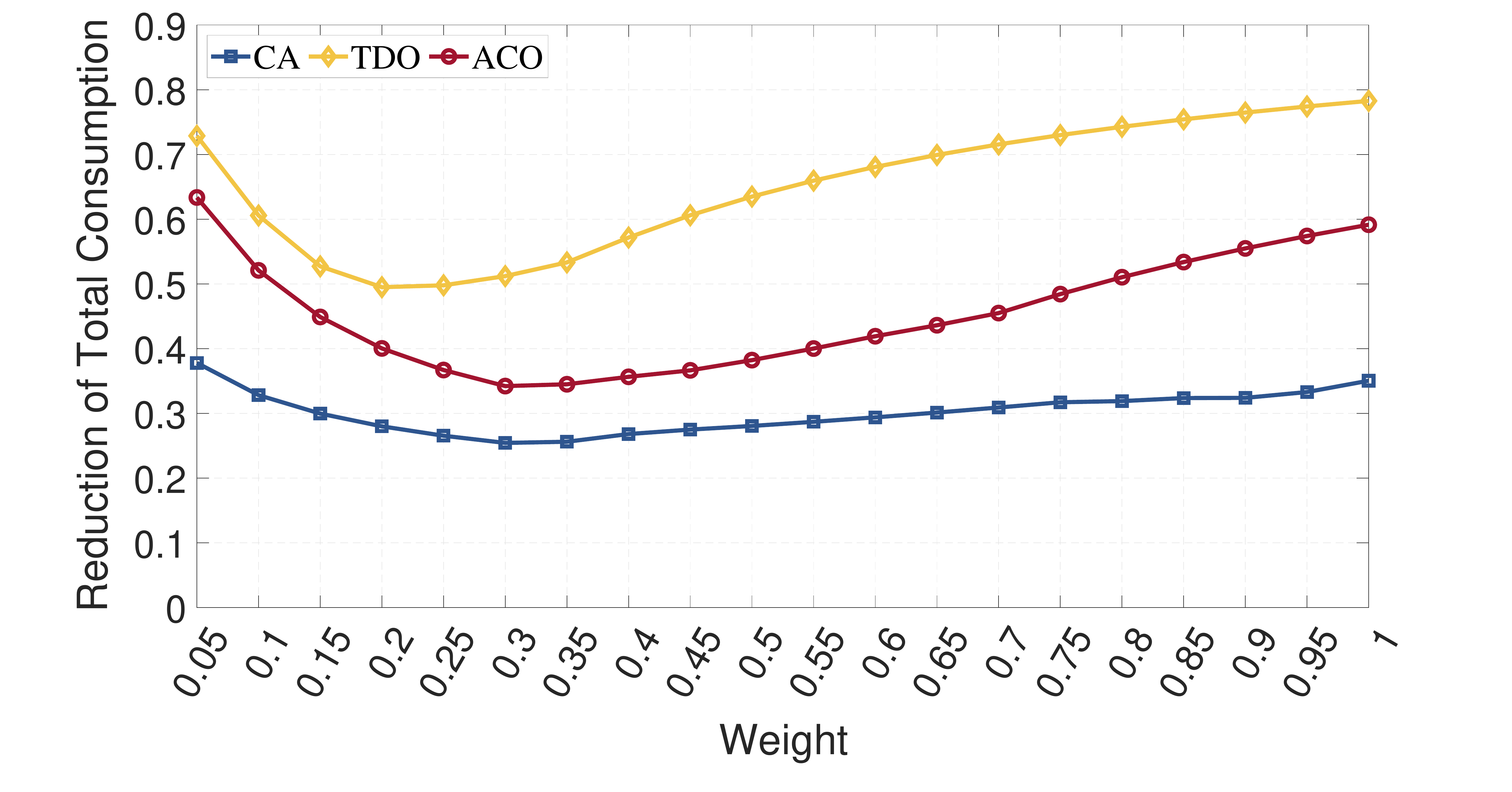}
    \caption{Relative reduction in $\mathcal{L}$ with different weights $\epsilon$}
    \label{fig:m result2}
    \vspace{-0.2cm}
\end{figure}

\subsubsection{Overall Performance Comparison}
Fig. \ref{fig:m result2} shows the performance comparison of CA, TDO and ACO with different values of the weighting factor $\epsilon$. The vertical axis in Fig. \ref{fig:m result2} is the relative reduction in $\mathcal{L}$ achieved by a certain policy $\Phi$ as compared to the baseline case RAN (($\mathcal{L}_{\text{RAN}}$-$\mathcal{L}_{\Phi}$)/$\mathcal{L}_{\text{RAN}}$, $\Phi \in \{$CA, TDO, ACO$\}$). For CA, TDO as well as RAN, the amounts of energy consumption are calculated based on the scenario where no obstacles are present, and the UAV moves along a straight-line trajectory. The results show that ACO incurs extra energy consumption over CA as the total travelling distance is longer for the UAV to avoid obstacles. However, ACO still achieves significant energy saving, up to 25.56\% compared to CA and 57.47\% compared to RAN. Therefore, the TDO is an effective complement to the obstacle avoidance mechanism in the trajectory planning phase that can be deployed to offset the extra energy consumption. The different trends of $ \mathcal{L}$ with changing weights are due to the scaling differences between energy and delay in the total consumption. 

\section{Conclusion}~\label{sec: conclusion}
In this paper, we proposed a framework to optimize the energy efficiency for the process of assigning computational tasks in a fog computing system, where a single UAV is available to assist in transmitting data and computing tasks. The proposed solution consists of two novel algorithms, namely a Particle Swarm Optimization based algorithm aiming at assigning tasks to different devices in the network, and an Ant Colony Optimization based algorithm to determine the optimal trajectory for the UAV. A greedy approach is then invoked to combine the outputs of the two algorithms to determine the location on the trajectory for the UAV to start transmission. 
We demonstrated that our proposed framework could considerably reduce the total power consumption compared to existing state-of-art methods while completing the same number of tasks under the same set of constraints. Furthermore, our resource allocation method converges faster and requires less computing power for the UAV compared to existing approaches. 

In the future, we plan to extend the model to account for the scenario where multiple UAVs are simultaneously deployed to assist computational offloading in fog computing. 

\section*{Acknowledgements}

This work is partly supported by Zhuhai Basic and Applied Basic Research Foundation Grant ZH22017003200018PWC, and partly supported by the Guangdong Provincial Key Laboratory of 
Interdisciplinary Research and Application for Data Science, BNU-HKBU United International 
College, Project code 2022B1212010006 and in part by Guangdong Higher Education Upgrading
Plan (2021-2025) UIC  R0400001-22.

\bibliographystyle{IEEEtran}
\bibliography{IEEEabrv, mybibliography} 
\end{document}